\documentclass[12pt]{amsart}
\usepackage{amsmath,amsthm}

\def\pmatrix{\left(\begin{matrix}}
\def\endpmatrix{\end{matrix}\right)}

\def\Sp{\operatorname{Sp}}

\def\R{{\mathbb R}}

\def\Z{{\mathbb Z}}
\def\F{{\mathbb F}_2}
\def\C{{\mathbb C}}
\def\N{{\mathbb N}}

\def\cal{\mathcal}

\def\SL{\operatorname{SL}}

\def\Coeff{\operatorname{Coeff}}

\def\tr{{\rm tr}}

\def\inv{^{-1}}
\def\<{\langle}
\def\>{\rangle}

\def\Half{{\mathcal H}}
\def\half{{\textstyle{\frac12}}}

\newcommand\smtwomat[4]{
{\bigl(
\genfrac{}{}{0pt}{1}{#1}{#3}\,
\genfrac{}{}{0pt}{1}{#2}{#4}
\bigr)}}

\def\2{\displaystyle}
\def\EE{{\cal E}}
\def\SS{{\cal S}_U}
\def\BB{{\cal B}}

\def\NN{{\mathbb N}}

\def\RR{{\mathbb R}}
\def\FF{{\mathbb F}}
\def\mymat#1#2#3#4{\left(\begin{matrix}#1&#2\\#3&#4\end{matrix}\right)}
\def\DU{\coprod}

\theoremstyle{plain}
\newtheorem{thm}{Theorem}
\newtheorem{lm}[thm]{Lemma}
\newtheorem{prop}[thm]{Proposition}
\newtheorem{cor}[thm]{Corollary}

\newtheorem{df}[thm]{Definition}

\theoremstyle{definition}
\begin{document}
\title[Binary forms]{ Binary forms 
and the hyperelliptic superstring Ansatz}

\author[C. Poor]{Cris Poor}
\address{Department of Mathematics, Fordham University, Bronx, NY 10458-5165, USA}
\email{poor@fordham.edu}

\author[D. Yuen]{David S. Yuen}
\address{Department of Mathematics and Computer Science, Lake Forest College, 555 N. Sheridan Rd., Lake Forest, IL 60045, USA}
\email{yuen@lakeforest.edu}



\begin{abstract}
We give a hyperelliptic formulation of the Ansatz of D'Hoker and Phong.  
We give an explicit family of binary invariants, one for each genus, 
that satisfies this hyperelliptic Ansatz.  
We also prove that this is the unique family 
of weight~eight binary forms over the theta group on the hyperelliptic locus
that satisfies this Ansatz.
Futhermore, we prove that
this solution may also be obtained 
by applying Thomae's map to multivalued Siegel modular forms of Grushevsky 
and making certain choices of roots.
\end{abstract}
\maketitle

\section{Introduction}

We formulate the Ansatz of D'Hoker and Phong for the ring of binary invariants, 
which can be viewed as a ring of modular forms on the moduli space of hyperelliptic Riemann 
surfaces.  We prove the existence and the uniqueness of the sequence 
of binary invariants $H_g$ satisfying the Ansatz.  Finally, we relate our work to a 
sequence of multivalued Siegel modular forms constructed by Grushevsky.  
When Thomae's formula is applied to Grushevsky's multivalued Siegel modular forms, each 
$H_g$ may be extracted as a certain branch.  

We first review the formulation of the Ansatz of D'Hoker and Phong on the Siegel 
upper half space $\Half_g$, where the description of the Witt map $\Psi^{*}_{ij}$ is simpler.  
The Ansatz has three parts.  
For each genus $g$, we seek Siegel modular forms of weight eight for the theta group, 
$\Xi^{(g)}[0] \in [\Gamma_g(1,2), 8]$, such that: $i)$~For all $g_1, g_2 \in \N$, 
$$
\Xi^{(g_1+g_2)}[0] \smtwomat{\Omega_1}{0}{0}{\Omega_2}= 
\Xi^{(g_1)}[0] (\Omega_1)\, \Xi^{(g_2)}[0](\Omega_2),   
$$
whenever $\Omega_i \in \Half_{g_i}$ are the period matrices of compact Riemann surfaces.  
We can rephrase this condition in terms of the Witt map, 
$\Psi^{*}_{i,j}: [\Gamma_{i+j}(1,2), 8]\to [\Gamma_i(1,2), 8] \otimes [\Gamma_j(1,2), 8]$, 
by saying  $\Psi^{*}_{g_1,g_2}\Xi^{(g_1+g_2)}[0] 
 =\Xi^{(g_1)}[0] \otimes \Xi^{(g_2)}[0] $ 
on the Jacobian locus.  
$ii)$~The trace of $\Xi^{(g)}[0] $ to level one, 
$\tr\left( \Xi^{(g_1+g_2)}[0]  \right) \in [\Gamma_g,8]$, vanishes on all 
$\Omega \in \Half_g$ that are period matrices of compact Riemann surfaces.  
$iii)$~The family of solutions to conditions~$i$ and $ii$ is uniquely determined by 
the genus one solution $\Xi^{(1)}[0]=\theta_0^4\,\eta^{12}$.  
This formulation of the Ansatz differs only slightly from the original by D'Hoker and Phong 
and its evolution may be traced in 
\cite{DHP1}, \cite{DHPa}, \cite{DHPb}, \cite{DHPc}, \cite{CDPvG}, \cite{CDvG2}, \cite{GR}, \cite{RSM} and \cite{OPSY}.  
The solutions for $g \le 2,3,4$~and~$5$ may be found in, for example, \cite{DHPa}, \cite{CDPvG},  \cite{GR} and \cite{OPSY}, respectively.  
Uniqueness is known for $g \le 4$.  It appears likely that the solution is also unique 
in $g=5$ and that for $g \ge 6$ no solutions exist.  

These mathematical questions owe their origin to the physics literature.  
We thank R. Salvati~Manni for introducing us to these ideas.  
The chiral superstring measure $d\nu[e]$ for a fixed theta characteristic~$e$ 
should take the form $d\nu[e]=f[e]^{(g)} d\mu$, where $d \mu$ is the Mumford measure 
and $f[e]^{(g)} $ is a weight eight Teichmuller modular form on the moduli space 
of curves with a fixed theta characteristic~$e$.  Condition~$i$ says that the measure 
should be the product measure on reducible curves.  Condition~$ii$ says that the 
traced level one measure, whose integral over moduli space gives the 
cosmological constant, vanishes pointwise.  These conditions are only required 
for period matrices of compact Riemann surfaces because the original interest is 
in Teichmuller modular forms on the moduli space of curves.  For $g \le 3$, 
period matrices are dense in $\Half_g$ but for $g \ge 4$ there is no a priori 
reason to expect that a solution $f[e]^{(g)} $ on Teichmuller space will analytically 
extend to all of $\Half_g$.  Thus it is remarkable that in $g=4$~and~$5$ the solutions 
$\Xi^{(g)}[e]$ exist as Siegel modular forms at all;  whereas the nonexistence of the 
$\Xi^{(g)}[e]$ for $g \ge 6$ would come as no surprise.  
The general existence of the Teichmuller forms  $f[e]^{(g)} $ remains open and has not 
even received a strict mathematical formulation--- a task best reserved for those who 
make significant progress.  Still, the above considerations have shown 
what the  $f[e]^{(g)} $ should be in $g \le 5$ and the existence of these 
$\Xi^{(g)}[0]$  is a remarkable vindication of the Ansatz of D'Hoker and Phong.  
For an entry into the physics literature see \cite{MorozovA}.  
For Teichmuller modular forms, see \cite{Ichikawa}\cite{IchikawaHE}.  

Another probe into the existence of the hypothetical $f[e]^{(g)} $ would be to restrict 
them to hyperelliptic curves, a special case that is always easier to study.  
If such a family exists on the moduli space of curves then it should also exist on the 
moduli space of hyperelliptic curves, although the uniqueness property might be lost.  
This idea is not new.  
In \cite{MorozovB},  A{.} Morozov studies the restriction of the $\Xi^{(g)}[e]$ 
to the hyperelliptic locus and recommends the general application of Thomae's formula 
to Grushevsky's multivalued Siegel modular form---  accomplished here in section~4.  
We give further vindication of the Ansatz of D'Hoker and Phong by 
formulating it for the moduli space of hyperelliptic curves and by proving that this 
formulation of the Ansatz is uniquely solvable.  
The form of the Witt map is more complicated in the hyperelliptic case but it can 
be found in Tsuyumine's work \cite{Tsu}.  The discussion of these broad topics ends with 
this Introduction but one can hope that having an explicit {\sl hyperelliptic approximation\/} 
to a chiral superstring measure for every genus will be of use.

The vector space of binary invariants of weight~$w$ in $r$~variables, 
$S_{w}(r)$, consists of those polynomials $f \in \C[a_1, \dots, a_r]$ satisfying 
\begin{align*}
\label{C1}
\forall\, \mymat ABCD\in\SL_2(\RR) &\text{ with } \gamma(z) =\frac{Az+B}{Cz+D}, \\
f(a_1,\ldots,a_r) &=f(\gamma(a_1),\ldots,\gamma(a_r))\prod_{i=1}^r(Ca_i+D)^w.
\end{align*}
From the matrix $\smtwomat{\sqrt{\lambda}}00{1/\sqrt{\lambda}}$ we see that each nontrivial 
$f \in S_w(r)$ is homogeneous of total degree $wr/2$ and from the matrix 
$\smtwomat{1}{\lambda}0{1}$ that each $f \in S_w(r)$ is a polynomial in the 
$a_i-a_j$.   We remark that any product of the $(a_i-a_j)$ where all of the $a_i$ 
occur exactly $w$~times is an element of $S_w(r)$.  For example, if we set 
$\Delta_T=\prod_{i,j \in T:\, i>j}(a_i-a_j)$ for $T \subseteq \{1, 2, \dots, r\}$, 
then $\Delta_{ \{1, 2, \dots, r\}}$ is an element of $S_{r-1}(r)$.  
The graded ring $S(r)=\oplus_{w=0}^{\infty} S_w(r)$ is integrally closed 
and finitely generated over~$\C$.  
We define the star map 
$*: \C[a_1, \dots, a_r] \to \C[a_1, \dots, a_{r-1}]$ by letting $f^*$ be the 
coefficient of the highest power of $a_r$ in $f$;  this makes $*$ a multiplicative map.  
Furthermore, $*$~is injective on $S_w(r)$.

In the context of binary invariants $S(2g+2)$, the theta group corresponds 
to a certain subgroup of the symmetric group $S_{2g+2}$.  
We also call the subgroup 
of permutations, ${\mathcal S}_U$,  which stabilizes the partition of $\{1,\dots,2g+2\}$ into even and odd 
elements, the {\it theta group\/} although any 
conjugate group would serve equally well.  The subspace of $S_w(2g+2)$ 
fixed elementwise by ${\mathcal S}_U$ is written $ S_w(2g+2)({\mathcal S}_U)$ .  
In the following section we will define a certain subspace 
 $\BB_g^k \subseteq S_{\frac12 k g}(2g+2)({\mathcal S}_U)$; 
 suffice it to say here that  $\BB_g^k$ is the largest subspace for which 
 applications of the Witt map 
 $W_{g_1, g_2}:  \BB_{g_1+g_2}^k \to  S_{\frac12 k g_1}(2g_1+2) \otimes  S_{\frac12 k g_2}(2g_2+2)$ 
 are defined.  Following  \cite{Tsu},  define the Witt map  
 for $f \in \BB_{g_1+g_2}^k $ by 
 \begin{align*}
&\left( (*\otimes *) W_{g_1, g_2}(f) \right)(a_1,\ldots,a_{2g_1+1},\alpha_1,\ldots,\alpha_{2g_2+1}) =  \\
&\text{Coeff}(t^{\frac12 k (2g_1g_2+g_1+g_2)}, f(a_1,\ldots,a_{2g_1+1},\alpha_1+t,\ldots,\alpha_{2g_2+1}+t)).
\end{align*}
We can now state a hyperelliptic Ansatz, modeled after that of D'Hoker and Phong.  
\smallskip

{\bf Hyperelliptic Ansatz\/} \newline 
We wish to find a sequence of binary invariants 
$H_g \in \BB_g^8$ such that
\begin{itemize}
\item i) 
$\forall g_1, g_2 \in \N: g=g_1+g_2, W_{g_1, g_2}(H_g)= H_{g_1} \otimes H_{g_2}$. 
\item ii) The symmetrization of  $H_g$ vanishes: 
$\sum_{\sigma \in S_{2g+2}} \sigma(H_g)=0$.  
\item iii) Any solution $H_g$ to i) and~{ii)} is uniquely determined by 
the genus one solution 
$H_1=\Delta_{ \{1,2,3,4 \} }(a_1-a_3)(a_2-a_4)$.  
\end{itemize}

A solution of i) and ii) for the original Anatz is taken to a solution of i) and ii) for the hyperelliptic Ansatz 
by Igusa's $\rho$-map.  The relevant commutative diagrams may be found in the 
final section.  In order to present the solution to this hyperelliptic Ansatz, 
we need the following definitions.  
\begin{df}
For a finite sequence of natural numbers $e=(e_1,\dots,e_r)$, 
define $\psi'_e = \prod_{i=1}^{r-1}(a_{e_i}-a_{e_{i+1}})$
and
$\psi_e=(a_{e_r}-a_{e_1})\psi'_e$.  \smallskip

Define 
$\EE_r=\{(e_1,\ldots,e_r)\in \N^r: \forall i, e_i\equiv i\mod 2 \text{ \rm and }\{e_1,\ldots,e_r\}= \{1,\dots,r\}\}$.
That is, $\EE_r$ consists of all  permutations of $(1,\ldots,r)$
that alternate odd and even, beginning with odd.
For $g \in \N$, let
$$
H_g = \frac{1}{2^g}  \frac{1}{g+1} \left(\Delta_{\{1,2,\dots,2g+2\}}\right)^2\sum_{e\in\EE_{2g+2}}\frac{-1}{\psi_e}.
$$
\end{df}

\begin{thm} 
The sequence $H_g\in \BB_g^8$ satisfies all three conditions of the hyperelliptic Ansatz. 
\end{thm}

\section{The Theorem}
This section requires the 
following additional notation.

\begin{itemize}
\item $B_g= \{1,2,\dots,2g+2\} $, $U_g=\{1,3,5,\ldots,2g+1\}$, $U'_g=B_g\backslash U_g$.
\item $\hat\EE_r=\{e\in\EE_r: e_r=r\}$.
\item For $e\in\EE_r$, define $e^*\in\EE_{r-1}$ to be the sequence 
obtained by deleting the last term in $e$.
Note that $e\mapsto e^*$ gives a natural isomorphism $\hat\EE_r\to\EE_{r-1}$.
\item Define $\SS = \{\sigma: \text{ $\sigma$  a permutation of } B_g: \sigma(U_g)=U_g \text{ or } U'_g\}$,
$\tilde\SS = \{\sigma: \text{ $\sigma$  a permutation of } B_g: \sigma(U_g)=U_g\}$.
\item For $f(a_1,\ldots,a_r)$ a polynomial and $\sigma$ a permuation,
define $\sigma(f)$ $ = f(a_{\sigma(1)},\ldots,a_{\sigma(r)})$.
\end{itemize}

For $T\subseteq B_g$, we let $T'=B_g\setminus T$ denote the complement of $T$ in $B_g$.  
When $\vert T \vert =g+1$, we note that $\Delta_{T}\Delta_{T'} \in S_g(2g+2)$.  
In fact, the ring $S^{(g)}(2g+2)= \oplus_{j=0}^{\infty} S_{gj}(2g+2)$ is the 
integral closure of the ring generated by the $\Delta_{T}\Delta_{T'}$ 
over all $T$ with $B_g= T \DU T'$ and 
$\vert T \vert =\vert T' \vert =g+1$, compare Igusa \cite{IgusaMFPI},  page 845, supplement I.  
For many purposes, this characterization of $S^{(g)}(2g+2)$ obviates the need to 
treat this ring abstractly.  
\begin{lm}
\label{C3}
A nontrivial $f \in S_w(r)$ has degree $w$ in each $a_i$.  The star map 
$*:S_w(r)\to \C[a_i-a_j; 1 \le i<j \le r-1]$ injects.  
\end{lm}
\begin{proof}
(Tsuyumine \cite{Tsu})  
Consider $\smtwomat{0}{1/\epsilon}{-\epsilon}{a_r} \in \SL_2(\C)$ for $\epsilon \ne 1$.  
We have $f(a_1,\dots,a_r)=$ 
$$
\left( a_r(1-\epsilon)  \right)^w  
f\left(\frac{1/\epsilon}{a_r-\epsilon a_1},\dots,
\frac{1/\epsilon}{a_{r}-\epsilon a_{r-1}},\frac{1/\epsilon}{a_r(1-\epsilon)}\right)
\prod_{i=1}^{r-1}(a_r-\epsilon a_{i})^w.  
$$
We let $\epsilon \to 1$ on both sides on this equation; 
the limit of the left hand side is the nontrivial polynomial~$f$.  
The limit of the right hand side does not exist if $\deg_{a_r} f >w$ and 
is zero if $\deg_{a_r} f <w$.   Thus $\deg_{a_r} f =w$ and the same holds for 
each variable~$a_i$.  The injectivity of the star map follows from taking the limit:  
$$
f(a_1,\dots,a_r)=
\left( \prod_{i=1}^{r-1}(a_r- a_{i})^w \right)  
f^*\left(\frac{1 }{a_r- a_1},\dots,
\frac{1 }{a_{r}-  a_{r-1}}\right)
.  
$$
To show that the polynomial $f^*$ lies in $\C[a_i-a_j; 1 \le i,j \le r]$, 
it suffices to check its invariance under translations: 
$f^*(a_1+\lambda, \dots, a_{r-1} +\lambda)=$
$\lim_{t\to\infty} t^{-w} f(a_1+\lambda, \dots, a_{r-1} +\lambda,t)=$ 
$\lim_{t\to\infty} t^{-w} f(a_1, \dots, a_{r-1} ,t-\lambda)=   f^*(a_1, \dots, a_{r-1} )$.   
\end{proof}

\begin{prop}
\label{C4}
For any $r\in\N$, $*\Delta_{\{1,2,\dots,r\}} = \Delta_{\{1,2,\dots,r-1\}}$.
For $T \subseteq B_g$ with $\vert T \vert =g+1$, 
$\left( \Delta_T \Delta_{T'} \right)^* = 
\begin{cases}
\Delta_{T \setminus\{2g+2\}}  \Delta_{T'}, \text{if $2g+2 \in T$,}  \\
\Delta_{T }  \Delta_{T'\setminus\{2g+2\}}, \text{if $2g+2 \not\in T$.} 
\end{cases}$

For $e\in\hat\EE_r$, we have
$*\psi_e = -\psi'_{e^*}$.
\end{prop}

\begin{proof}
The proof is straightforward.
\end{proof}

For $r=r_1+r_2$ and $j \in \N$, we follow Tsuyumine by defining a map
\begin{align*}
T_{r_1,r_2}^{(j)}:\, &\C[a_1,\dots,a_{r}]  \to \C[a_1,\dots,a_{r_1}] \otimes \C[\alpha_1,\dots,\alpha_{r_2}]  \\
f \mapsto &
\text{Coeff}(t^j, f(a_1,\ldots,a_{r_1},\alpha_1+t,\ldots,\alpha_{r_2}+t))
\end{align*}
and a valuation 
\begin{align*}
\nu_{r_1,r_2}:\,  &\C[a_1,\dots,a_{r}]  \to \Z  \\
f \mapsto &
\deg_t f(a_1,\ldots,a_{r_1},\alpha_1+t,\ldots,\alpha_{r_2}+t).  
\end{align*}

\begin{df}
Define a valuation subring $S(2g+2)_0$ by 
\begin{align*}
&S_w(2g+2)_0=
\{f \in S_w(2g+2):     \\
&\forall g_1, g_2 \in \N: g_1+g_2=g, 
\nu_{2g_1+1,2g_2+1}( f) \le \frac{w}{g}(2g_1g_2 +g_1 +g_2) \}.  
\end{align*}
\end{df}

\begin{lm}
\label{C5}
Let $g_1, g_2 \in \N$ with $g=g_1+g_2$.  For $T \subseteq B_g$ define 
$\pi_1 T= \{x: x\in T \text{ and } 1 \le  x \le 2g_1+1\}  \subseteq B_{g_1}$ and 
$\pi_2 T= \{x-(2g_1+1): x\in T \text{ and } 2g_1+2 \le  x \le 2g+2\} \subseteq B_{g_2}$.  
For $T \subseteq B_g$ with $\vert T \vert =  g+1$, we have 
$\Delta_T\Delta_{T'}\in S_g(2g+2)_0$ and 
$T_{2g_1+1,2g_2+1}^{(2g_1g_2+g_1+g_2)}( \Delta_T\Delta_{T'} )=$
$$
\begin{cases}
0,\quad \text{\rm if } \vert \pi_1 T \vert \not\in \{g_1+1, g_1\}, \\
\left( \Delta_{\pi_1 T}  \Delta_{(\pi_1 T)'} \right)^* \otimes \left( \Delta_{(\pi_2 T')'}  \Delta_{\pi_2 T'} \right)^*, 
\quad\text{\rm if } \vert \pi_1 T \vert = g_1+1, \\
\left( \Delta_{\pi_1 T'}  \Delta_{(\pi_1 T')'} \right)^* \otimes \left( \Delta_{(\pi_2 T)'}  \Delta_{\pi_2 T} \right)^*, 
\quad\text{\rm if } \vert \pi_1 T \vert = g_1.  
\end{cases}
$$
\end{lm}
\begin{proof}
Let $m=  \vert \pi_1 T \vert $ and $n=  \vert \pi_1 T' \vert $.  
Then $m+n=2g_1+1$ and $\nu_{2g_1+1, 2g_2+1}\left(  \Delta_T\Delta_{T'} \right)=m(g+1-m)
+n(g+1-n)= (2g_1+1)(g+1)-(m^2+n^2)$.  
In this case, $m,n\in \Z_{\ge 0}$ and $m+n$ is odd, 
so the minimum of $m^2+n^2$ occurs at $\{m,n\}=\{g_1+1, g_1\}$.  
Therefore, 
 $\nu_{2g_1+1, 2g_2+1}\left(  \Delta_T\Delta_{T'} \right) \le  (2g_1+1)(g+1)-((g_1+1)^2+g_1^2)
 = 2g_1g_2+g_1+g_2$ and $\Delta_T\Delta_{T'}\in S_g(2g+2)_0$.  
 
 To find the coefficient of $t^{ 2g_1g_2+g_1+g_2 }$ in the cases of equality 
 we may assume $ \vert \pi_1 T \vert = g_1+1$ and $ \vert \pi_1 T '\vert = g_1$;  
 the other case follows by swapping $T$ and $T'$.   We have 
 \begin{align*} 
 & \Delta_T(a_1,\dots,a_{2g_1+1}, \alpha_1+t,\dots,\alpha_{2g_2+1}+t)=  \\
 & \prod_{i,j \in \pi_1 T: i > j} (a_i-a_j) 
 \prod_{i,j \in \pi_2 T: i > j} (\alpha_i-\alpha_j) 
 \prod_{i \in \pi_2 T, j \in \pi_1 T} (\alpha_i +t -a_j) 
 \end{align*}
 and similarly for $\Delta_{T'}$ so that 
  \begin{align*} 
 & T_{2g_1+1,2g_2+1}^{(2g_1g_2+g_1+g_2)}( \Delta_T\Delta_{T'} )(a_1,\dots,a_{2g_1+1}, \alpha_1,\dots,\alpha_{2g_2+1})= \\
 & \Delta_{\pi_1 T}(a_1,\dots,a_{2g_1+1}) \Delta_{\pi_2 T}(\alpha_1,\dots,\alpha_{2g_2+1}) \\
 & \Delta_{\pi_1 T'}(a_1,\dots,a_{2g_1+1}) \Delta_{\pi_2 T'}(\alpha_1,\dots,\alpha_{2g_2+1}) = \\
& \Delta_{\pi_1 T}(a_1,\dots,a_{2g_1+1}) 
\Delta_{(\pi_1 T)' \setminus \{ 2g_1+2\} } (a_1,\dots,a_{2g_1+1}) \\
& \Delta_{(\pi_2 T')'  \setminus \{ 2g_2+2\} } (\alpha_1,\dots,\alpha_{2g_2+1})
\Delta_{\pi_2 T'  } (\alpha_1,\dots,\alpha_{2g_2+1}) .  
 \end{align*}
 Thus, $T_{2g_1+1,2g_2+1}^{(2g_1g_2+g_1+g_2)}( \Delta_T\Delta_{T'} )=$ 
 $\left( \Delta_{\pi_1 T} \Delta_{(\pi_1 T)'} \right)^* \otimes 
 \left( \Delta_{(\pi_2 T')'} \Delta_{\pi_2 T' }\right)^* $ 
 upon comparison with Proposition~\ref{C4}.
\end{proof}

\begin{cor}
\label{C6}
The map 
$$
T_{2g_1+1,2g_2+1}: S^{(g)}(2g+2)_0 \to 
\oplus_{j=0}^{\infty} \,\,S_{g_1 j}(2g_1+2)^* \otimes S_{g_2 j}(2g_2+2)^*
$$
defined by 
$$
T_{2g_1+1,2g_2+1}^{ (j( 2g_1g_2+g_1+g_2 )) }: S_{g j}(2g+2)_0 \to 
 S_{g_1 j}(2g_1+2)^* \otimes S_{g_2 j}(2g_2+2)^*
$$
is a homomorphism of graded rings.  
\end{cor}
\begin{proof}
We need to check that the codomain is as stated.  
The previous Lemma~\ref{C5} shows this for the ring generated by the $ \Delta_T\Delta_{T'}$; 
thus it holds for any subring of the integral closure where 
$T_{2g_1+1,2g_2+1}$ is multiplicative.  We know that 
$T_{2g_1+1,2g_2+1}$ is multiplicative on $S^{(g)}(2g+2)_0$ 
by the valuation condition defining $S^{(g)}(2g+2)_0$.  
\end{proof}
Since the star map is injective, the Witt map $W_{g_1,g_2}$ 
is well-defined by the following:
\begin{df}
\label{C7}
Let $g_1,g_2\in \N$ with $g_1+g_2=g$.  
The graded ring homomorphism 
$ W_{g_1, g_2}: S^{(g)}(2g+2)_0 \to \oplus_{j=0}^{\infty} 
\,S_{g_1 j}(2g_1+2)  \otimes S_{g_2 j}(2g_2+2)  $
is defined on $S_{gj}(2g+2)_0$ by 
$(* \otimes *)\circ W_{g_1,g_2} = T_{2g_1+1,2g_2+1}^{(j (2g_1g_2+g_1+g_2))}$.  
\end{df}

Intuitively, the $T$ map pulls apart a hyperelliptic surface and the star map opens up
a hyperelliptic surface at a branch point; so the Witt map pulls apart a hyperelliptic
surface into two pieces and then closes up the individual pieces.

\begin{prop}\label{prop1a}
Given $g=g_1+g_2$, with $g_1,g_2\in\NN$,
we have
\begin{itemize}
\item
$\nu_{2g_1+1,2g_2+1}\Delta_{B_g}
=(2g_1+1)(2g_2+1)$ 
and
$T_{2g_1+1,2g_2+1}^{(4g_1g_2+2g_1+2g_2+1)}\Delta_{B_g}
=
\Delta_{ \{ 1,\dots,{2g_1+1} \} } \otimes 
\Delta_{ \{ 1,\dots,{2g_2+1} \} } 
$
\item
$\nu_{2g_1+1,2g_2+1} \Delta_{U_g}
=(g_1+1)g_2$ and
$\nu_{2g_1+1,2g_2+1} \Delta_{U'_g} 
=g_1(g_2+1)$
\item
If $e\in\EE_{2g+2}$ such that $\{e_1,\ldots,e_{2g_1+1}\}={ \{ 1,\dots,{2g_1+1} \} }$,
define
$e_L=(e_1,\ldots,e_{2g_1+1})$
and
$e_R=(e_{2g_1+2}-(2g_1+1),\ldots,e_{2g+2}-(2g_1+1))$.
Then
$T_{2g_1+1,2g_2+1}^{(2)}\psi_e
=
-\psi'_{e_L} \otimes 
\psi'_{e_R}.$
\end{itemize}
\end{prop}
\begin{proof}
The proof is straightforward.
\end{proof}

\begin{prop}\label{prop1b}
Let $g=g_1+g_2$, with $g_1,g_2,j \in\NN$.  
Denote $\delta=j(2g_1g_2+g_1+g_2)$ and let $h\in S_{gj}(2g+2)_0$.
Supose $h=f_1f_2$ with both $f_1,f_2$ polynomials.
Suppose $\nu_{2g_1+1,2g_2+1} f_1 
=\delta_1$. 
Then
$T_{2g_1+1,2g_2+1}^{(\delta)}h=
T_{2g_1+1,2g_2+1}^{(\delta_1)}f_1\cdot
T_{2g_1+1,2g_2+1}^{(\delta-\delta_1)}f_2$
\end{prop}

\begin{proof}
The fact that $h$ satisfies the valuation property implies
that
$\nu_{2g_1+1,2g_2+1} f_2 
\le \delta-\delta_1$.
The result then follows easily.
\end{proof}

\begin{df}
\label{C10}
For subgroups $G$ of the symmetric group $S_{2g+2}$, define 
$$
S_w(2g+2)_{0}(G)= 
\{ f\in S_w(2g+2)_{0}: \forall \sigma \in G, \sigma(f)=f \}.  
$$
For even $k$, define $\BB_g^k= S_{\frac12 k g}(2g+2)_{0}(\SS)$.  
\end{df}
This space of binary invariants, $\BB_g^k$,  is analogous to 
$[\Gamma_g(1,2),k]$, the Siegel modular forms of degree~$g$ and 
weight~$k$ for the theta group.   It remains to define the concept of a cusp form on $\BB_g^k$.  
For $1 \le m,n \le r$, 
define $\bar\Phi_{mn}$ on a polynomial $f\in \C[a_1,\ldots,a_r]$ 
by $\bar\Phi_{mn}f = f$ with $a_m=0=a_n$.
For $m\ne n$, define $\Phi_{mn}: S_{gj}(2g+2) \to  \C[a_i;  1 \le i \le 2g+2, i \ne m,n]$ by
$\Phi_{mn} f = \bar\Phi_{mn}f/(\prod_{\ell\ne m,n}a_\ell)^{-j}$
and reindexing the variables if necessary.

\begin{lm}
\label{C11}
Let $m,n,j\in \N$.  We have 
$ \Phi_{mn}: S_{gj}(2g+2) \to S_{(g-1)j}(2g)$.  
\end{lm}
\begin{proof}
Consider  $\Delta_T\Delta_{T'} \in S_g(2g+2)$.  
If $\{m,n\} \subseteq T$ or $\{m,n\} \subseteq T'$ then $ \Phi_{mn}(\Delta_T\Delta_{T'} )=0$.  
Otherwise we may relabel so that $m \in T$ and $n \in T'$ and then 
$ \Phi_{mn}(\Delta_T\Delta_{T'} )= \pm \Delta_{T\setminus \{ m \} }\Delta_{T' \setminus \{n \}}$, 
which  is indeed in $S_{g-1}(2(g-1)+2)$ after potential reindexing.  
Because the $ \Phi_{mn}$ map is multiplicative, the codomain is shown to be as stated 
by taking integral closure.  
\end{proof}

\begin{df}
\label{C12}
An element $f \in S^{(g)}(2g+2)$ is called a cusp form if 
$ \Phi_{mn}( f )=0$ for all distinct $m,n \in B_g$.   
\end{df}

\begin{thm}{\label{main} {\rm (Main Theorem)\/}}
For each $g\in\N$, define
\begin{equation}
H_g =\frac{1}{2^g}\frac{1}{g+1} \left(\Delta_{B_g}\right)^2  \sum_{e\in\EE_{2g+2}}\frac{-1}{\psi_e}.
\end{equation}
Then the following conditions hold:  
\begin{itemize}
\item $H_1=\Delta_{B_1}\Delta_{U_1}\Delta_{U'_1}
$.
\item $H_g \in \BB^8_g$.
\item For all $g_1,g_2\in\N$ with $g_1+g_2=g$,
we have that
$$W_{g_1,g_2}H_g = H_{g_1}\otimes H_{g_2}$$
\item $\sum_{\sigma\in S_{2g+2}}\sigma(H_g) = 0$.
\item $\Phi_{ij}H_g=0$ for all $i\ne j$.
\end{itemize} 
\end{thm}

\begin{proof}
It is straightforward to check that
$H_1=\Delta_{B_1}\Delta_{U_1}\Delta_{U'_1}$.  
The last two conditions are also easily checked:  
Consider the polynomial 
$G= \frac1{2^g} \frac{1}{g+1} \Delta_{B_g}\sum_{e\in\EE_g}\frac{-1}{\psi_e}$, so that 
$H_g=\Delta_{B_g}G$.
For any two $i\ne j$, 
$\Phi_{ij}H_g = \Phi_{ij}\Delta_{B_g}\cdot
\Phi_{ij}(G) 
= 0 \cdot \Phi_{ij}(G)
=0.$  
Next, we prove that $\sum_{\sigma\in S_{2g+2}}\sigma(H_g) $ is trivial.  
We have
$$\sum_{\sigma\in S_{2g+2}}\sigma(H_g) = 
\frac{1}{2^g} 
\frac{1}{g+1}\Delta_{B_g}^2\sum_{\sigma\in S_{2g+2}}\sum_{e\in\EE_g}\frac{-1}{\sigma(\psi_e)}
$$
because $\Delta_{B_g}^2$ is invariant under all $\sigma\in S_{2g+2}$.
Define a polynomial by $\tilde G=\frac{1}{2^g} \Delta_{B_g}\frac{1}{g+1}\sum_{\sigma\in S_{2g+2}}
\sum_{e\in\EE_g}\frac{-1}{\sigma(\psi_e)}$.  
Since $\Delta_{B_g}$ is alternating and since $\Delta_{B_g}\tilde G$ is invariant
under $S_{2g+2}$,
then $\tilde G$ must be alternating.
This implies that $\tilde G$ is a multiple of $\Delta_{B_g}$.
But $\deg \tilde G < \deg \Delta_{B_g}$ forces $\tilde G=0$.

Next, we show that $H_g\in\BB^8_g$.
First,  it is clear from the construction that
$H_g\in S_{w}(2g+2)$ where $w=4g$
because it is a sum whose terms are products of the form $(a_i-a_j)$
where $i$, $j$ are of opposite parity such that in the product each $a_i$
appears exactly $w$ times.
Second, $H_g$ is invariant under $\SS$ because
$$\sum_{e\in\EE_{2g+2}}\frac{1}{\psi_e}=
\sum_{\sigma\in\tilde\SS}\frac{1}{\sigma(\psi_{e_0})}$$ 
for any particular $e_0\in\SS$,
and because applying any $\tau\in \SS$ we have
$\tau(\psi_{e_0})=\psi_{e_1}$ for some $e_1\in\EE_{2g+2}$.
Third,  the valuation property will be evident when 
we find the image of the Witt map.  

From the definition of the Witt map, we need to prove that
$$
\text{Coeff}(t^{8g_1g_2+4g_1+4g_2}, 
H_g(a_1,\ldots,a_{2g_1+1},\alpha_1+t,\ldots,\alpha_{2g_2+1}+t)
),
$$ 
is equal to
$(*H_{g_1})(a_1,\ldots,a_{2g_1+1})\cdot
 (*H_{g_2})(\alpha_1,\ldots,\alpha_{2g_2+1})$.
Note from Proposition~\ref{prop1a} that the maximal power of $t$ in the expansion of the factor 
$\Delta_{B_g}^2(a_1,\ldots,a_{2g_1+1},\alpha_1+t,\ldots,\alpha_{2g_2+1}+t)$
is $t^{2(2g_1+1)(2g_2+1)}$
and that its coefficient is 
$\Delta_{\{ 1,\dots,{2g_1+1} \} }^2(a_1,\ldots,a_{2g_1+1})\cdot
\Delta_{ \{ 1,\dots,{2g_2+1} \} }^2(\alpha_1,\ldots,\alpha_{2g_2+1})$.
We claim that $\nu_{2g_1+1,2g_2+1}(\psi_e) \ge 2$.  
For simplicity, let $C=\{1,\dots,a_{2g_1+1}\}$ and $D=B_g\backslash C$.
Then $\nu_{2g_1+1,2g_2+1}(\psi_e)=$
$\deg_t(\psi_e(a_1,\ldots,a_{2g_1+1},\alpha_1+t,\ldots,\alpha_{2g_2+1}+t))$
is the number of transitions between the sets $C$ and $D$
in the sequence
$e_1,\ldots,e_r,e_1$.
This number is clearly at least 2
and is exactly 2 if and only if all the numbers in $C$ are together and all the numbers in $D$
are together (where we have to view the sequence with wrap-around);
call the set of such $e$ the set $\cal F$.
In particular, we just proved that $H_g$ satisfies
the valuation condition 
$\nu_{2g_1+1,2g_2+1}(H_g) \le 2(2g_1+1)(2g_2+1) -2=4(2g_1g_2+g_1+g_2)$ 
so that $H_g \in S_{4g}(2g+2)_0(\SS)=\BB_g^8$.  

Since  
$ t^{8g_1g_2+4g_1+4g_2} = t^{2(2g_1+1)(2g_2+1)} / t^2$,
if $
\nu_{2g_1+1,2g_2+1}\psi_e >2 
$ for an~$e$, 
then $\text{Coeff}(t^{8g_1g_2+4g_1+4g_2},
\Delta_{B_g}^2/\psi_e(a_1,\ldots,a_{2g_1+1},\alpha_1+t,\ldots,\alpha_{2g_2+1}+t)
=0$.
Thus we have that
$$
T_{2g_1+1,2g_2+1}^{(8g_1g_2+4g_1+4g_2)}H_g =
T_{2g_1+1,2g_2+1}^{(8g_1g_2+4g_1+4g_2+2)}\left(\Delta_{B_g}^2\right)\cdot
\frac{2^{-g}}{g+1}\sum_{e\in\cal F}\frac{-1}{T_{2g_1+1,2g_2+1}^{(2)}\psi_e}
,
$$
Since $\psi_e$ is unchanged when $e$ is cyclically rotated,
we may rotate $e$ so that the set $C$ comes first and then the set $D$.
To this end,  
let
$$\tilde{{\cal F}}=\{e\in\EE_{2g+2}:
\{e_1,\ldots,e_{2g_1+1}\}=\{1,\dots,a_{2g_1+1}\} \}.$$
Since there are $g+1$ ways to cycle $e$ from an element of $\tilde{{\cal F}}$
to an element of ${\cal F}$,
we can replace the sum over ${\cal F}$ by a sum over $\tilde{{\cal F}}$ 
and an overall factor of $g+1$:  
$$
T_{2g_1+1,2g_2+1}^{(8g_1g_2+4g_1+4g_2)}H_g =
T_{2g_1+1,2g_2+1}^{(8g_1g_2+4g_1+4g_2+2)}  \frac{1}{2^g} \left(\Delta_{B_g}^2\right)
\sum_{e\in\tilde{\cal F}}\frac{-1}{T_{2g_1+1,2g_2+1}^{(2)}\psi_e}.
$$
Then by Proposition~\ref{prop1a}   
\begin{align*}
T_{2g_1+1,2g_2+1}H_g & {=}
\Delta_{\{1,\dots,{2g_1+1}\} }^2(a_1,\ldots,a_{2g_1+1}) 
\Delta_{\{1,\dots,{2g_2+1}\} }^2(\alpha_1,\ldots,\alpha_{2g_2+1}) \\
 \frac{1}{2^{g_1+g_2}} 
&\sum_{e_L\in\EE_{2g_1+1}}\sum_{e_R\in\EE_{2g_2+1}}\frac{-1\cdot-1}{\psi'_{e_L}(a_1,\ldots,a_{2g_1+1})
\psi'_{e_R}(\alpha_1,\ldots,\alpha_{2g_2+1}}).
\end{align*}
because we can view
each $e\in\tilde{\cal F}$
as the concatenation of two pieces $e_L$ and $e_R$;
that is given an $e\in\tilde{\cal F}$,
we have corresponding
$e_L=(e_1,\ldots,e_{2g_1+1})$
and
$e_R=(e_{2g_1+2}-(2g_1+1),\ldots,e_{2g+2}-(2g_1+1))$.

On the other hand,
\begin{align*}
*H_{g_1} = & \frac{1}{2^{g_1}} 
\left(*\Delta_{B_{g_1}}\right)^2 \frac{1}{{g_1}+1}\sum_{e\in\EE_{2g_1+2}}\frac{-1}{*\psi_e}\\
=&  \frac{1}{2^{g_1}}  \left(\Delta_{ \{1,\dots,{2g_1+1} \} }\right)^2 \sum_{e\in\hat\EE_{2g_1+2}}\frac{-1}{*\psi_e}\\
\end{align*}
$$
=  \frac{1}{2^{g_1}}  \left(\Delta_{ \{1,\dots,{2g_1+1}\} }\right)^2 \sum_{e\in\hat\EE_{2g_1+2}}\frac{-1}{-\psi'_{e^*}}
=  \frac{1}{2^{g_1}} \left(\Delta_{ \{1,\dots,{2g_1+1} \} }\right)^2 \sum_{e\in\EE_{2g_1+1}}\frac{1}{\psi'_{e}}
$$
and similarly for $*H_{g_2}$.

Now it is easy to see that 
$T_{2g_1+1,2g_2+1}H_g = *H_{g_1} \otimes *H_{g_2}.$
\end{proof}

We thank R. Salvati  Manni for bringing the following consequence to our attention: 
As A{.} Morozov points out in \cite{MorozovB}, the fact that $\Delta_{B_g} $ divides $H_g$ 
implies that, for variables $x$ and $y$ and $P_T(x)=\prod_{i \in T} (x-a_i)$, 
$$
\sum_{\sigma \in S_{2g+2}} \sigma \left( 
(P_{U}(x) P_{U'}(y)- P_{U}(y) P_{U'}(x) ) H_g \right) =0.   
$$
The reason for this is that the complete  symmetrization must be divisible by $(x-y)\Delta_{B_g}^2$.  
Along with $\sum \sigma(H_g)=0$, 
this identity is equivalent to the 
non-renormalization of the 2 and 3-point functions.  

\section{Uniqueness}
We now prove some propositions aimed at proving the uniqueness of the family $H_g$.

\begin{prop}\label{prop11}
For any $r\in\N$,
let $f\in S_w(r)$.
Then
\begin{equation*}
(T_{r-3,3}^{(3w)}f)(a_1,\ldots,a_{r-3}, \alpha_1,\alpha_2,\alpha_3) = \text{Coeff}(t^{3w}, f(a_1,\ldots,a_{r-3},t,t,t)).
\end{equation*}
Furthermore,  $f(a_1,\ldots,a_{r-3},u,u,u)=$ 
\begin{equation*}
\prod_{i=1}^{r-3}(u-a_i)^w
\cdot
(T_{r-3,3}^{(3w)}f)(\frac1{u-a_1},\ldots,\frac1{u-a_{r-3}} , \alpha_1,\alpha_2,\alpha_3  ).
\end{equation*}
\end{prop}

\begin{proof}
Since each variable in $f$ occurs to degree $w$,
it is clear that 
\begin{multline}
\text{Coeff}(t^{3w}, f(a_1,\ldots,a_{r-3},t+\alpha_1,t+\alpha_2,t+\alpha_3)=  \\
\text{Coeff}(t^{3w}, f(a_1,\ldots,a_{r-3},t,t,t))
\end{multline}
and that this expression is really independent of $\alpha_1,\alpha_2,\alpha_3$.  
Therefore
$(T_{r-3,3}^{(3w)}f)(a_1,\ldots,a_{r-3},\alpha_1,\alpha_2,\alpha_3) = \text{Coeff}(t^{3w}, f(a_1,\ldots,a_{r-3},t,t,t))$.

Fix $0< \epsilon <1 $. Let $u$ be a variable.
Let $\gamma(z)=\frac{1/\epsilon}{u-\epsilon z}$.
Then $f\in S_w(r)$ implies that 
\begin{align*}
f(a_1,\ldots,a_{r-3},u,u,u) = &
\prod_{i=1}^{r-3}(u-\epsilon a_i)^w\cdot\\
(u-\epsilon u)^{3w} 
&
f(\tfrac{1/\epsilon}{u-\epsilon a_1},\ldots,
\tfrac{1/\epsilon}{u-\epsilon a_{r-3}},\tfrac{1/\epsilon}{(1-\epsilon) u},
\tfrac{1/\epsilon}{(1-\epsilon) u},\tfrac{1/\epsilon}{(1-\epsilon) u}
).  
\end{align*}
Let us expand $f(a_1,\ldots,a_{r-3},t,t,t)$ in powers of $t$ as
\begin{multline}
f(a_1,\ldots,a_{r-3},t,t,t) = \\
(T_{r-3,3}^{(3w)}f)(a_1,\ldots,a_{r-3},\alpha_1,\alpha_2,\alpha_3) t^{3w}
+ G(a_1,\ldots,a_{r-3},t),
\end{multline}
where $\deg_tG<3w$.
Then
\begin{align*}
f(a_1  &,\ldots,a_{r-3},u,u,u) \\
=&\prod_{i=1}^{r-3}(u-\epsilon a_i)^w\cdot
(
(T_{r-3,3}^{(3w)}f)(\tfrac{1/\epsilon}{u-\epsilon a_1},\ldots,\
\tfrac{1/\epsilon}{u-\epsilon a_{r-3}}, \alpha_1,\alpha_2,\alpha_3) (1/\epsilon)^{3w}\\
&+(1-\epsilon)^{3w} G(\tfrac{1/\epsilon}{u-\epsilon a_1},\ldots,
\tfrac{1/\epsilon}{u-\epsilon a_{r-3}},\tfrac{1/\epsilon}{(1-\epsilon) u})
)\\
=&\prod_{i=1}^{r-3}(u-\epsilon a_i)^w\cdot
(
(T_{r-3,3}^{(3w)}f)(\tfrac{1/\epsilon}{u-\epsilon a_1},\ldots,\
\tfrac{1/\epsilon}{u-\epsilon a_{r-3}},\alpha_1,\alpha_2,\alpha_3) (1/\epsilon)^{3w}\\
&+(1-\epsilon)^{3w} (\text{terms where $(1-\epsilon)^\beta$ occurs with $\beta>-3w$})
)
\end{align*}
Taking the limit as $\epsilon\to1$ gives the desired result.
\end{proof}

\begin{prop}\label{prop12}
Let $f\in S_{w}(2g+2)(\SS)$.
If $f(a_1,\ldots,a_{2g+2})=0$ whenever
$a_i=a_j=a_k$ with distinct $i,j,k$ not all of the same parity, 
then either $f=0$ or $\deg f\ge g(g+1)$.
\end{prop}

\begin{proof}
Assume $f\ne0$.
For each integer $0\le j\le g+1$,
define a polynomial $h_j$ by
\begin{multline*}
h_j(x_1,y_1,\ldots,x_j,y_j,b_1,\ldots,b_{g+1-j})=\\
f(x_1,y_1,\ldots,x_j,y_j,b_1,b_1,\ldots,b_{g+1-j},b_{g+1-j}).
\end{multline*}
Note
$\deg f\ge \deg h_j$ for each $j$.
Note that $h_{g+1}=f$,
so $h_{g+1}\ne 0$.
Then let $m$ be the minimum such that $h_m\ne0$.
Since $f$ is invariant under $\SS$,
then $h_m$ is invariant under swapping within the $x_i$ or within the $y_i$,
and $h_m$ is invariant under swapping within the $b_i$.
Note $h_m=0$ whenever $b_i=b_j$ with $i\ne j$,
Thus
$$
h_m=\prod_{0\le i<j\le g+1-m}(b_i-b_j)\cdot
 k(x_1,y_1,\ldots,x_j,y_j,b_1,\ldots,b_{g+1-j}),
$$
for some polynomial $k$.
Then $k$ would be alternating under swapping within the $b_i$
which implies that $k$ is a multiple of each $(b_i-b_j)$.
Thus
$$
h_m=\prod_{0\le i<j\le g+1-m}(b_i-b_j)^2\cdot
 k_2(x_1,y_1,\ldots,x_j,y_j,b_1,\ldots,b_{g+1-j}),
$$
for some polynomial $k_2$.
Now, also $h_m=0$ whenever any $b_i=x_j$ or $b_i=y_j$.
Thus
\begin{multline*}
h_m=\prod_{0\le i<j\le g+1-m}(b_i-b_j)^2\cdot
\prod_{i,j}(b_i-x_j)(b_i-y_j)\cdot\\
 k_3(x_1,y_1,\ldots,x_j,y_j,b_1,\ldots,b_{g+1-j}),
\end{multline*}
for some polynomial $k_3$.
Then the homogeneous degree is 
$$\deg h_m\ge (g+1-m)(g-m)+2(g+1-m)m=(g+1-m)(g+m).$$
If $m=0$, then this says $\deg h_0\ge (g+1)g$
and $\deg f\ge (g+1)g$ follows.
If $m>0$, then $h_{m-1}=0$.
This says that $h_m=0$ if $x_m=y_m$.
Thus $h_m$ is a multiple of $(x_m-y_m)$, and
so $h_m$ is a multiple $(x_i-y_j)$ for all $i,j$.
Thus
\begin{multline*}
h_m=\prod_{0\le i<j\le g+1-m}(b_i-b_j)^2\cdot
\prod_{i,j}(b_i-x_j)(b_i-y_j)\cdot\\
\prod_{i,j}(x_i-y_j)(b_i-y_j)\cdot
 k_4(x_1,y_1,\ldots,x_j,y_j,b_1,\ldots,b_{g+1-j}),
\end{multline*}
for some polynomial $k_4$.
Then
$$\deg h_m\ge (g+1-m)(g+m)+m^2=(g+1)g+m.$$
Then $\deg f\ge (g+1)g+m$ and the proposition is proved.
\end{proof}

\begin{cor}\label{cor12}
For any $g\in\N$ with $g\ge2$, 
let $f\in S_{g-1}(2g+2)(\SS)$.
If
$T_{2g-1,3}^{(3g-3)}f= 0$,
then  $f=0$.
\end{cor}

\begin{proof}
Suppose we have an $f\in S_{g-1}(2g+2)(\SS)$ with 
$T_{2g-1,3}^{(3g-3)}f= 0$.
Proposition \ref{prop11} with $w=g-1$ implies
that $f(a_1,\ldots,a_{2g-1},u,u,u)=0$.
By symmetry under $\SS$,
this implies $f(a_1,\ldots,a_{2g+2})=0$ whenever
three of the $a_i$ are equal with not all three indices of the same parity.
By Proposition \ref{prop12}, we have
either $f=0$ or
$\deg f\ge (g+1)g$.
But if $f\ne0$, then $f\in S_{g-1}(2g+2)$ implies that
$\deg f= \frac12(g-1)(2g+2)=g^2-1<(g+1)g$, a contradiction.
Hence $f=0$.
\end{proof}

\begin{prop}\label{prop7}
Any cusp form $h\in\BB^k_g$ must be of the form
$$h=\Delta_{B_g}\Delta_{U_g}\Delta_{U'_g}f,$$
where $f\in S_{\frac12kg -3g-1}(2g+2)$.
\end{prop}

\begin{proof}
Since $\Phi_{ij}h=0$ for any $i\ne j$,
then $h=0$ whenever $a_i=a_j$.
This forces $(a_i-a_j)$ to be a divisor of $h$.
Thus $h=\Delta_{B_g} h_2$ for some polynomial $h_2$.
Since $h$ is invariant under $\SS$ and $\Delta_{B_g}$
is alternating under $\SS$,
then $h_2$ must be alternating under $\SS$,
which means that $h_2$ changes sign whenever $a_i$ and
$a_j$ are swapped with $i,j$ of the same parity.  
This implies $h_2=0$ whenever $a_i=a_j$ with $i,j$ of the same parity.
So $h_2$ must be a multiple of $\Delta_{U_g}$ and $\Delta_{U'_g}$.
Thus $h=\Delta_{B_g}\Delta_{U_g}\Delta_{U'_g}f$
with $f$ a polynomial.
Since $h$ has  weight $\frac12 kg$
and $\Delta_{B_g}\Delta_{U_g}\Delta_{U'_g}$
has  weight $3g+1$,
then $f$ has the asserted weight.  
\end{proof}

\begin{lm}\label{lemma7}
Let $f$ be a  cusp form in $S^{(g)}(2g+2)$.  
For $i,j<2g+2$, we have $\bar\Phi_{ij} (*f)=0$ .
\end{lm}

\begin{proof}
Let $f \in S_{g \ell}(2g+2)$.  
We have 
\begin{align*}
\bar\Phi_{ij}(*f)=&\Phi_{ij} \Coeff(t^{g\ell},f(a_1,\dots,a_{2g+1}, t))
\\
=&\text{$\Coeff(t^{g\ell},f(a_1,\dots,a_{2g+1}, t))$ with $a_i=a_j=0$ }
\\=&\text{$\Coeff(t^{g\ell},f(a_1,\dots,a_{2g+1}, t)$ with $a_i=a_j=0$ })
\\=&0.  
\end{align*}
\end{proof}

\begin{prop}\label{phiprop}
Let $f\in\BB_g^k$ be a binary invariant with respect to the theta group.
Suppose $W_{g-1,1}f=h_2\otimes h_1$.
Then $*\Phi_{2g+1,2g+2}f =  (-1)^{k/2} (*h_2) \bar\Phi_{2,3}(*h_1)\alpha_1^{-k/2}$
 and 
$*\Phi_{2g,2g+2}f =   (*h_2) \bar\Phi_{1,3}(*h_1)(-\alpha_2)^{-k/2}$.
In particular, if $W_{g-1,1}f=h_2\otimes h_1$
where $h_1$ is a cusp form, then $f$ is a cusp form.
\end{prop}

\begin{proof}
Just write out 
\begin{align*}
(*h_2)(a_1,&\ldots,a_{2g-1})\,(*h_1)(\alpha_1,\alpha_2,\alpha_3)
= 
T_{2g-1,3}f
\\
= &
\Coeff(t^{\frac12k(3g-2)},f(a_1,\ldots,a_{2g-1},\alpha_1+t,\alpha_2+t,\alpha_3+t))
\end{align*}
Then
\begin{align*}
(*h_2)\,\bar\Phi_{23}(*h_1)
= &
\Coeff(t^{\frac12k(3g-2)},f(a_1,\ldots,a_{2g-1},\alpha_1+t,t,t))
\\
= &
\Coeff(t^{\frac12k(3g-2)},f(a_1-t,\ldots,a_{2g-1}-t,\alpha_1,0,0))
\end{align*}
On the other hand,
\begin{align*}
*\Phi_{2g+1,2g+2}f =&
*\left(
\frac{f(a_1,\ldots,a_{2g},0,0)}{(a_1\cdots a_{2g})^{k/2}}
\right)
\\
=&
\Coeff(t^{k(g-1)/2},
\frac{f(a_1,\ldots,a_{2g-1},t+a_{2g},0,0)}{(a_1\cdots a_{2g-1})^{k/2} (t+a_{2g})^{k/2}})
\\
=&
\Coeff(t^{k(g-1)/2},
\frac{f(a_1-t,\ldots,a_{2g-1}-t,a_{2g},0,0)}{((a_1-t)\cdots (a_{2g-1}-t))^{k/2} (a_{2g})^{k/2}})
\end{align*}
where we used the fact that $*\Phi_{2g+1,2g+2}f$ is invariant under translations
in the last equality using Lemmas~\ref{C3} and~\ref{C11}.    
 Since the highest term in~$t$ in the denominator
is $(-1)^{(2g-1)k/2} a_{2g}^{k/2} t^{(2g-1)k/2}$,
then
\begin{align*}
*\Phi_{2g+1,2g+2}f  &=
(-1)^{(2g-1)k/2}  (a_{2g})^{-k/2}\cdot
\\
&\Coeff(t^{k(g-1)/2+(2g-1)k/2},
f(a_1-t,\ldots,a_{2g-1}-t,a_{2g},0,0)
) 
\\
=(-1&)^{k/2}  (a_{2g})^{-k/2}\cdot
\\
&\Coeff(t^{k(3g-2)/2},
f(a_1-t,\ldots,a_{2g-1}-t,a_{2g},0,0)
) 
\end{align*}
This proves
$*\Phi_{2g+1,2g+2}f = (-1)^{k/2} (*h_2) \bar\Phi_{2,3}(*h_1)\alpha_1^{-k/2}$
and similarly
$*\Phi_{2g,2g+2}f =  (-1)^{k/2} (*h_2) \bar\Phi_{1,3}(*h_1)\alpha_2^{-k/2}$.

Now suppose  that  $W_{g-1,1}f=h_2\otimes h_1$
where  $h_1$ is a cusp form.
Then by Lemma \ref{lemma7},
${\bar \Phi}_{2,3}(*h_1)=0$ and so $*\Phi_{2g+1,2g+2}f=0$
and thus $\Phi_{2g+1,2g+2}f=0$.
Similarly $\Phi_{2g,2g+2}f=0$.
The invariance of~$f$ under $\SS$ implies
$\Phi_{i,j}f=0$ for all $i\ne j$ and so $f$ is a cusp form.
\end{proof}

\begin{prop}\label{prop9}
The Witt map
$W_{g-1,1}$ is injective on $\BB_g^8$.  
\end{prop}

\begin{proof}
Let $h\in\BB^8_g$ and suppose $W_{g-1,1}h=0$.
Then $T^{(12g-8)}_{2g-1,3}h=* (W_{g-1,1}h)=0$.  
By Proposition~\ref{phiprop} we deduce that $h$ is a cusp form.  
By Proposition \ref{prop7},
we know that any cusp form $h\in \BB_g^8$
must be of the form
$h=\Delta_{B_g}\Delta_{U_g}\Delta_{U'_g}f$
where $f\in S_{g-1}(2g+2)$.
From Proposition~\ref{prop1a},
we know $\nu_{g-1,1}(\Delta_{B_g}\Delta_{U_g}\Delta_{U'_g})
=
(2g-1)3+(g)1+(g-1)2=9g-5.$
Then Proposition \ref{prop1b}
says that
\begin{equation}
T^{(12g-8)}_{2g-1,3} h
= T^{(9g-5)}_{2g-1,3}(\Delta_{B_{g}}\Delta_{U_{g}}\Delta_{U'_{g}})
\cdot
T_{2g-1,3}^{(3g-3)}f .
\end{equation}
Since $T^{(12g-8)}_{2g-1,3} h=0$ and
$T^{(9g-5)}_{2g-1,3}(\Delta_{B_{g}}\Delta_{U_{g}}\Delta_{U'_{g}})\ne0$,
then $T_{2g-1,3}^{(3g-3)}f=0$.
Since $f\in S_{g-1}(2g+2)$, then Corollary \ref{cor12}
implies $f=0$.
So $h=0$, completing the proof.
\end{proof}

\begin{thm}\label{unique}
The family $H_g$ as given in Theorem \ref{main}
is the unique family that satisfies the
first three conditions  stated in that Theorem.
\end{thm}

\begin{proof}
Suppose by way of contradiction there is another family $K_g$ that 
satisfies the first three conditions of Theorem \ref{main}.
Let $g_0$ be the smallest index such that
$K_{g_0}\ne H_{g_0}$.
Then $g_0\ge 2$ by the first condition.  
Use the third condition to check 
\begin{align*}
W_{g_0-1,1}(K_{g_0}-H_{g_0})
&=W_{g_0-1,1}K_{g_0}-W_{g_0-1,1}H_{g_0}\\
&=K_{g_0-1}\otimes K_{1}-H_{g_0-1}\otimes H_{1}\\
&=0.
\end{align*} 
Since $W_{g_0-1,1}$ is injective on $\BB_{g_0}^8$ 
by Proposition \ref{prop9}, 
and $K_{g_0}-H_{g_0} \in \BB_{g_0}^8$ by the second condition, 
then $K_{g_0}-H_{g_0}=0$,
which is a contradiction.
\end{proof}

\section{Remarks on Grushevsky's Construction}

In \cite{GR}, Grushevsky gave a uniform construction of  Siegel modular cusp forms
that satisfied the Ansatz in genera $g=1,2,3,4$:
\begin{equation}
\Xi^{(g)}[0]=
\frac1{2^g}\left(\sum_{i=0}^g(-1)^i 2^{\frac12i(i-1)}G^{(g)}_{i,2^{4-i}}
\right)
\end{equation}
where
$$
G^{(g)}_{i,r}=\sum_{V\subseteq\FF_2^{2g}}\left(\prod_{\zeta\in V}
\theta[\zeta]\right)^r$$
and where the sum is over isotropic subspaces $V$ of dimension $i$.  

Since $\Xi^{(g)}[0]$ is multivalued for $g > 4$, 
it is natural to ask whether some branch is single valued on the Jacobian locus.  
In \cite{GSM2}, S. Grushevsky and R. Salvati~Manni showed that, if single valued, 
$\tr (\Xi^{(g)}[0])$ is a multiple of $J^{(g)}$, 
the difference of the theta series of the two classes of even unimodular rank~{16} 
lattices.  For $1 \le g \le 3$, $J^{(g)}$ is trivial whereas $J^{(4)}$ is the Schottky form 
defining the Jacobian locus, see \cite{Ch}.  The long open problem of whether 
$J^{(g)}$ vanishes on the Jacobian locus for $g>4$ was resolved negatively in  \cite{GSM2}; 
thus $\Xi^{(g)}[0]$ stops solving the Ansatz for $g > 4$.  
However, it is known \cite{Poor} that $J^{(g)}$ always vanishes on the hyperelliptic locus 
and we will show that $\rho_g\left(  \Xi^{(g)}[0] \right)$ does have a branch that 
solves the hyperelliptic Ansatz.  Thus, the intricate pattern discovered by 
Grushevsky in the construction of $\Xi^{(g)}[0]$ properly belongs to the 
hyperelliptic locus even though the same pattern happens to define  
a Siegel modular form for $g \le 4$.  We need some definitions and lemmas.  

We refer to \cite{Fr} and \cite{Mumford} for standard theory on Siegel modular forms.  
The action of $M=\smtwomat{A}BCD \in \Sp_g(\R)$ on $\Omega \in \Half_g$ is 
$M<\Omega>=(A\Omega+B)(C\Omega+D)\inv$.  Let $\Gamma \subseteq \Gamma_g=\Sp_g(\Z)$ 
be a subgroup of finite index.  The vector space of Siegel moduar forms of degree~$g$, 
weight~$k$ and character $\chi$, written $[\Gamma, \chi, k]$ is the set of holomorphic 
functions $f:\Half_g\to\C$, bounded at the cusps for $g=1$, that satisfy 
$f|_k M= \chi(M) f$ for all $M \in \Gamma$ where 
$(f|_k M)(\Omega)= \det(C\Omega+D)^{-k}f( M<\Omega> )$.  
We define the graded rings:  $M^{(k_0)}\left(\Gamma,\chi\right)= \sum_{j=0}^{\infty} [\Gamma, \chi^j, jk_0]$.  

The set $\{ S \subseteq B_g: \vert S \vert \text{ even } \}$ is a group under the 
symmetric difference $\oplus$.  The quotient group 
$\{ S \subseteq B_g: \vert S \vert \text{ even } \}/ \{ \emptyset, B_g \}$ treats each 
$S$ as equivalent to its complement~$S'$.  In fact, we have an explicit isomorphism 
\begin{align}\label{iso}
\eta: ( \{ S \subseteq B_g: \vert S \vert \text{ even } \} / \{ \emptyset, B_g \},\oplus) 
&\to ( \F^{2g}, +)
\end{align}
given by $\eta_S=\sum_{i\in S}\eta_i$ and 
\begin{alignat*}{2}
&\eta_1=\begin{matrix}1&0&\cdots&0\\0&0&\cdots&0\end{matrix}\,,  &
&
\eta_2=\begin{matrix}1&0&\cdots&0\\1&0&\cdots&0\end{matrix}\,,
\\
&\vdots
   &  &
\vdots
\\
&\eta_{2i-1}=\begin{matrix}0&\cdots&0&1&0&\cdots&0\\1&\cdots&1&0&0&\cdots&0\end{matrix}\,,  \qquad&
&
\eta_{2i}=\begin{matrix}0&\cdots&0&1&0&\cdots&0\\1&\cdots&1&1&0&\cdots&0\end{matrix}\,,
\\
\end{alignat*}
\begin{alignat*}{2}
& \vdots
&  &
\vdots
\\
&\eta_{2g-1}=\begin{matrix}0&0&\cdots&0&1\\1&1&\cdots&1&0&\end{matrix}, \qquad
&
&\quad \eta_{2g}=\begin{matrix}0&0&\cdots&0&1\\1&1&\cdots&1&1&\end{matrix},
\\
&\eta_{2g+1}=\begin{matrix}0&0&\cdots&0&0&\\1&1&\cdots&1&1&&\end{matrix},  \qquad
&
&\quad \eta_{2g+2}=\begin{matrix}0&0&\cdots&0&0&\\0&0&\cdots&0&0&&\end{matrix}.
\end{alignat*}

We treat the elements $\zeta \in \F^{2g}$ as theta characteristics; i.e., 
we consider the action $\zeta \mapsto M\cdot \zeta$ of $\Sp_g(\F)$ on $\F^{2g}$ given by
\begin{equation}
\label{affine}
\theta[\zeta] \vert_{1/2} M \in (\text{eighth roots of unity}) \,\theta[M \cdot \zeta].
\end{equation}
Explicitly, we have 
$$
\smtwomat{A}BCD \cdot \zeta = \smtwomat{A}BC{D}' \zeta 
+\begin{pmatrix} (A'C)_0 \\ (B'D)_0 \end{pmatrix}
$$
where $(X)_0$ denotes the vector formed from the diagonal of $X$ 
and $\zeta$ is treated as a column vector.  
In general this action is affine and is linear precisely when we have $M \in \Gamma_g(1,2)(\F)$; 
this can be taken as the definition of  $\Gamma_g(1,2)$.  
Frobenius found a complete set of 
invariants for this action, \cite{I}, page~{212}.   
For $\zeta=\begin{bmatrix} a \\ b \end{bmatrix}$, $\zeta_1,\zeta_2,\zeta_3\in \F^{2g}$, we put

\begin{align*}
&e_*(\zeta)=(-1)^{a\cdot b},                     \\
e(\zeta_1,\zeta_2,\zeta_3 )  &=
e_*(\zeta_1)e_*(\zeta_2)e_*(\zeta_3)
e_*(\zeta_1+\zeta_2+\zeta_3).
\end{align*}

The Theorem of Frobenius can be stated as follows:
\begin{thm} 
\label{F1}
Let $\Sp_g(\F)$ act on theta characteristics in $ \F^{2g}$ as in equation~\ref{affine}.  
Two sequences,  
$(\zeta_1,\dots,\zeta_m)$ and $(\xi_1,\dots,\xi_m)$, 
are in the same $\Sp_g(\F)$-orbit 
if and only if sending 
$\zeta_i \mapsto \xi_i$ preserves 
\begin{itemize}
\item all linear relations with an even number of summands, 
\item all $e_*$ values and 
\item all $e$ values.  
\end{itemize}
\end{thm}

Given any permutation~$\sigma$ of $B_g$, 
we can induce a linear map 
${\bar \sigma}: \F^{2g} \to  \F^{2g} $ by $\eta_S \mapsto \eta_{\sigma(S)}$ for 
 $S \subseteq B_g$ with $ \vert S \vert $ even.  
 This action is  induced by an element $M \in \Sp_g(\F)$ 
 if and only if $\sigma$ preserves $e_*$, in view of the linearity of ${\bar \sigma}$.   
 One can check, or see \cite{PoorHE}, page 824, that 
 $e_*(\eta_S)= (-1)^{ \frac12 (g+1- \vert S \oplus U \vert )  }$, 
 so that for $\sigma \in \SS$ there exists an $M \in \Sp_g(\F)$ 
 such that $\eta_{\sigma(S)}=M \cdot \eta_S$.  
 This $M$ is uniquely determined because the $\eta_{\{i,j\}}$ span $\F^{2g}$ 
 and  we have $M \in \Gamma_g(1,2)(\F)$ because ${\bar \sigma}$ is linear.  
 We will have use for a certain character on $\Gamma_g(1,2)$.  
 Define $\kappa$ by 
 $\theta[0]|M=\kappa\theta[0]$.  
 Then $\kappa^4$ gives a real character of $\Gamma_g(1,2)(\F)$, 
 or of $\Gamma_g(1,2)$.  From Igusa \cite{I}, page 182, we know that 
 $\kappa^4$ is given by 
$\smtwomat{A}BCD \mapsto (-1)^{\tr(D-I_g)}$.  

We now connect the traditional marking of a hyperelliptic curve with Igusa's 
$\rho$-homomorphism.  Let $ W \subseteq \C^{2g+2}$ be the quasiprojective 
variety of points with distinct coordinates.  There is a morphism 
$h: W \to \Gamma_g(2) \backslash \Half_g$ that sends 
$a=(a_1,\dots,a_{2g+2}) \in W$ to the $\Gamma_g(2)$-class 
of the period matrix $\Omega(a)$ for the traditional marking \cite{Mumford} of a hyperelliptic curve 
$y^2=\prod_{i=1}^{2g+2} (x-a_i)$.  
The $\rho_g$ map follows Thomae's formula, given below, 
and for all $f,g \in [\Gamma_g(2), k]$ with $\rho_g(g)\ne 0$ we have 
the important property, \cite{Tsu}, 
page 777.  
\begin{equation}
\label{Z3}
\dfrac{\rho_g(f)}{\rho_g(g)}=\dfrac{f\circ h}{g \circ h} .
\end{equation}

\begin{lm}
\label{Z1}
Let $\sigma \in \SS$ and $M \in \Gamma_g(1,2)$ with 
$\eta_{\sigma(S)} = M \cdot \eta_S$ for all $S \subseteq B_g$ with $\vert  S \vert$ even.  
For all $a \in W$, we have $h(a^{\sigma})=M\langle h(a) \rangle$.  
\end{lm}
\begin{proof}
Let ${\mathcal C\/}$ be the Riemann surface given by the hyperelliptic curve 
$y^2=\prod_{i=1}^{2g+2} (x-a_i)$.  
Let $w: {\mathcal C\/} \to \operatorname{Jac}({\mathcal C\/})= \C^g / 
\Lambda(\Omega(a))$ be the Abel-Jacobi map, where 
$\Lambda(\Omega)= \Z^g + \Omega \Z^g$.  
In the traditional marking of a hyperelliptic curve ${\mathcal C\/}$ we have 
$w( (a_i,0) )= \frac12 \left( \Omega(a),  I \right) \eta_i$, see \cite{PoorHE}, page 824,  
and the Lemma follows from this as we explain.  

The points $a, a^{\sigma} \in W$ both define ${\mathcal C\/}$  but the traditional markings,  
see page~{3{.}76} of \cite{Mumford},  will differ.  
Let $\begin{pmatrix} B \\ A \end{pmatrix}$  be the standard homology basis corresponding to 
$a \in W$ and $\Omega(a)$ the period matrix computed from this basis.  
Similarly,  let $\begin{pmatrix} {\tilde B} \\  {\tilde A} \end{pmatrix}$ correspond to 
$a^{\sigma} $ so that 
$\begin{pmatrix} {\tilde B} \\  {\tilde A} \end{pmatrix}
=\smtwomat{\alpha}{\beta}{\gamma}{\delta} \begin{pmatrix} B \\ A \end{pmatrix}$ for 
some $\smtwomat{\alpha}{\beta}{\gamma}{\delta} \in \Sp_g(\Z)$ and we have 
$\Omega( a^{\sigma} )= \smtwomat{\alpha}{\beta}{\gamma}{\delta}\langle \Omega(a) \rangle = 
(\alpha \Omega(a) + \beta) (\gamma \Omega(a) + \delta)\inv$.  
The Abel-Jacobi maps 
$w :  \operatorname{Div}^0\left( {\mathcal C\/} \right) \to \C^g / \Lambda( \Omega(a))$, 
${\tilde w} :  \operatorname{Div}^0\left( {\mathcal C\/} \right) \to \C^g / \Lambda( \Omega(a^{\sigma})) $,  
are related by 
${\tilde w} = ( \Omega(a) \gamma' + \delta')\inv w$.  
Thus we have
\begin{align*} 
 ( \Omega(a) \gamma' + \delta')\inv   &w( (a_{\sigma(i)},0)-(a_{\sigma(j)},0) )  \equiv 
{\tilde  w}( (a_i^{\sigma},0)-(a_{j}^{\sigma},0) )   \\
&\equiv 
\frac12 \left( \Omega(a^{\sigma}), I \right) \eta_{ \{i,j\} } \mod \Lambda( \Omega(a^{\sigma}) ),  
\end{align*}
\begin{align*} 
&\frac12 \left( \Omega(a ), I \right) M' \eta_{ \{i,j\} } \equiv 
\frac12 \left( \Omega(a ), I \right) M \cdot \eta_{ \{i,j\} } \equiv 
\frac12 \left( \Omega(a ), I \right)  \eta_{ \{ \sigma(i), \sigma(j) \} } \equiv \\
&w( (a_{\sigma(i)},0)-(a_{\sigma(j)},0) )  \equiv 
\frac12 \left((\Omega(a) \gamma' + \delta'  ) \Omega(a^{\sigma} ), 
\Omega(a) \gamma' + \delta'  \right)   \eta_{ \{i,j\} } \equiv  \\\
&\frac12 \left( \Omega(a ), I \right) \smtwomat{\alpha}{\beta}{\gamma}{\delta}' \eta_{ \{i,j\} } 
\mod  \Lambda( \Omega(a ) ).  
\end{align*}
Thus we have $ M \equiv  \smtwomat{\alpha}{\beta}{\gamma}{\delta} \mod 2$ and, 
along with $\Omega( a^{\sigma} )= \smtwomat{\alpha}{\beta}{\gamma}{\delta}\langle \Omega(a) \rangle $, 
this implies 
$\Gamma_g(2) \Omega(a^{\sigma}) = \Gamma_g(2) M\langle  \Omega(a ) \rangle$.  
Since $\Gamma_g(2)$ is normal in $\Gamma_g(1,2)$, this is  
$h(a^{\sigma} )=   M\langle  h(a ) \rangle $.  
\end{proof}

\begin{df}
Define a subset $S\subseteq\{1,\ldots,2g+2\}$
to be balanced if $S$ contains an equal number of even numbers and odd numbers,
and unbalanced otherwise.
\end{df}

Here is Thomae's formula:  
If $T\subset B_g$ with $|T|=g+1$, we have
\begin{equation}
\rho:\, \theta[\eta_{U_g\oplus T}]^4
\mapsto
(-1)^{\lfloor\frac{g+1}2\rfloor}\cdot
\prod_{i<j\in T}(a_i-a_j)
\prod_{i<j\in T'}(a_i-a_j) .  
\end{equation}
If $\zeta$ cannot be put into the form $\eta_{U_g\oplus T}$
with $|T|=g+1$,
then $\theta[\zeta]\mapsto0$.
It is simple to see that that $S=U_g\oplus T$ with $|T|=g+1$
if and only if $S$ is balanced.
Therefore, when $S$ is unbalanced,
then we have $\rho(\theta[\eta_S])=0$.
By the isomorphism (\ref{iso}),
we say that a theta characteristic is {balanced} if
it can be written as $\eta_S$ with $S$ balanced. 

We can now give the commutative diagrams that show that a solution of the Ansatz of D'Hoker and Phong 
goes to a solution of the hyperelliptic Ansatz under Igusa's 
$\rho$-homomorphism;  
 these commutative diagrams are deduced from those in Tsuyumine \cite{Tsu}.

\begin{prop}
\label{T1}
Let $g_1,g_2,g\in\N$ such that $g_1+g_2=g$.  
There exists an $M_{g_1,g_2}\in \Gamma_g(1,2)$ such that 
the following diagram commutes.  
\begin{align*}
M^{(2)}\left( \Gamma_g(2) \right) 
& \overset{{\hskip0.5in \rho_g \circ \vert M_{g_1,g_2} \hskip0.5in}}{\longrightarrow} S^{(g)}(2g+2)_{ 0}   \\
{\Psi_{g_1,g_2}^{*}} \downarrow & \hskip1.8in
 \downarrow{W_{g_1,g_2}} \\
M^{(2)}\left( \Gamma_{g_1}(2) \right)  {\otimes} M^{(2)}  &\left( \Gamma_{g_2}(2) \right)  
\overset{{\rho_{g_1} \otimes \rho_{g_2} }}{\longrightarrow}  
S^{(g_1)}(2g_1{+}2) {\otimes} S^{(g_2)}(2g_2{+}2)
\end{align*}
Let $\kappa^4$ be the character of $\Gamma_g(1,2)$ given by 
$\smtwomat{A}BCD \mapsto (-1)^{\tr(D-I_g)}$.  
 $$
\rho_g: M^{(2)}\left( \Gamma_g(1,2), \kappa^4 \right) 
\to \BB_g   .
$$
\end{prop}
\begin{proof}
The commutative diagram actually given in \cite{Tsu}, page 786, 
and there only in the case of $(g_1,g_2)=(g-1,1)$, is 
\begin{align*}
M^{(2)}\left( \Gamma_g(2) \right) 
& \overset{{\hskip0.7in \rho_g \hskip0.7in}}{\longrightarrow} S^{(g)}(2g+2)_{0}   \\
{\Psi_{g_1,g_2}^{*}} \downarrow & \hskip1.8in
 \downarrow{T_{2g_1+1,2g_2+1}} \\
M^{(2)}\left( \Gamma_{g_1}(2) \right)  {\otimes} M^{(2)}  &\left( \Gamma_{g_2}(2) \right)  
\overset{(*\otimes *){(\rho_{g_1} \otimes \rho_{g_2}) }}{\longrightarrow}  
S^{(g_1)}(2g_1{+}2)^* {\otimes} S^{(g_2)}(2g_2{+}2)^*
\end{align*}
However, the proof of the general case is the same.  
In this article, we fix $\rho_g$ to be the map induced by the traditional 
marking of a hyperelliptic curve;  therefore the commutative diagram from Tsuyumine 
holds for $ \rho_g \circ \vert M_{g_1,g_2} $ for some 
$M_{g_1,g_2}\in \Sp_g(\Z)$.  Applying this commutative diagram to 
$\theta[0]^8$ we see that $M_{g_1,g_2}\in \Gamma_g(1,2)$.  

In order to show that $\rho_g$ sends 
$M^{(2)}\left( \Gamma_g(1,2), \kappa^4 \right) $ to $\BB_g$, 
let $f \in [ \Gamma_g(1,2), \kappa^{2k} , k] $ and consider 
$\sigma \rho_g(f)$ for $ \sigma \in \SS$.  Notice first that 
$$
\sigma\left( \rho_g\left( \theta[0]^{2k} \right) \right) = 
\sigma\left( \epsilon_g \Delta_{U} \Delta_{U'} \right)^{k/2}= 
\left( \epsilon_g \Delta_{U} \Delta_{U'} \right)^{k/2}= 
\rho_g\left( \theta[0]^{2k} \right).  
$$
Therefore, using equation~\ref{Z3} and Lemma~\ref{Z1}, we have 
\begin{align*}
&\dfrac{ (\sigma( \rho_g(f)))(a) }{ (\sigma( \rho_g( \theta[0]^{2k})))(a) }= 
\dfrac{ f( h( a^{\sigma}) )}{\theta[0]^{2k}( h(a^{\sigma} ) )}=
\dfrac{ f(M \cdot h( a ) )}{\theta[0]^{2k}( M \cdot h(a  ) )}=  \\ 
&\dfrac{ (f\vert M)(  h( a ) )}{(\theta[0]^{2k}\vert M)(   h(a  ) )}=
\dfrac{ f(  h( a ) )}{ \theta[0]^{2k} (   h(a  ) )}=
\dfrac{ ( \rho_g(f))(a) }{  ( \rho_g( \theta[0]^{2k}))(a) }.  
\end{align*}
Thus, $\sigma \rho_g(f)=\rho_g(f)$.  
\end{proof}

\begin{lm}\label{balancelemma}
Let $V$ be an isotropic subspace with all balanced elements.
Then there exists a partitioning of $\{1,\ldots,2g+2\}$
into balanced subsets of 2 elements each (call
them $u_1,\ldots,u_{g+1}$),
and a subspace $H\subseteq \F^{g+1}$ with $\dim H=\dim V$
such that
$$
V=\{\eta_{S_h}: h\in H, \text{ where }S_h=\bigcup_{i: h_i\ne0} u_i \}
$$
\end{lm}

\begin{proof}
Using the isomorphism (\ref{iso}),
we view $V$ as a set of balanced subsets $\{S_j\}$.
Given that the symmetric difference $S_{j_1}\oplus S_{j_2}$
is balanced by hypothesis, then
the intersection $S_{j_1}\cap S_{j_2}$ is also balanced.
More generally, we can prove that the intersection of
any number of these balanced subsets will be balanced.
Then  the Venn Diagram of intersections of all the $S_j$
will give a partitioning of $\{1,\ldots,2g+2\}$
into balanced subsets.
We can make a finer partition into balanced subsets
of 2 elements each such that each $S_j$ is a union
of subcollection of this partition.
The result follows.
\end{proof}

\begin{lm}\label{thelemma}
Fix genus $g$ and fix an $d\in\NN$ with $0\le d\le g$.
For any subspace $V$ of dimension $d$,
there exists a polynomial $Q_V\in\Z[a_1,\dots,a_{2g+2}]$ in the $a_i$ such that
we have
\begin{equation*}
\rho(\prod_{\zeta\in V}\theta[\zeta]^8)= Q_V^{2^d}.
\end{equation*}
Furthermore, $Q_V^2$ is unique.  
\end{lm}

\begin{proof}
Note that if  $V$ is not isotropic or if $V$ contains
any unbalanced theta characteristics,
then $Q_V=0$ suffices.
So assume $V$ is isotropic and contains only balanced elements.
Let $H\subseteq \F^{g+1}$ and a partitioning of $\{1,\ldots,2g+2\}$
into balanced subsets $u_1,\ldots,u_{g+1}$
 of 2 elements each,
as in Lemma \ref{balancelemma}.
Note $|H|=2^d$.
Now, we have
$$\rho(\prod_{\zeta\in V}\theta[\zeta]^4)=\pm
\prod_{i>j}(a_i-a_j)^{r_{ij}}
$$
where $r_{ij}$ are exponents that we will calculate.  
There are two cases: $i,j$ of the same or different parity.

Case: $i,j$ are of the same parity.  Let $i\in u_a$ and $j\in u_b$
with $a\ne b$.
Then $\rho(\theta[\eta_{S_h}]^4)$ contains a factor of
$(a_i-a_j)$ if and only if $i,j$ are both in or both not in
$S_h\oplus U$ which happens (because they are of the same parity)
 if and only if they are both in or not in
$S_h$ which is if and only if $h_a=h_b$.
Since $H$ is a vector subspace of $\F^{g+1}$, then
the number of $h\in H$ for which $h_a=h_b$ is
either $|H|$ or $\half |H|$
because the set of such is the kernel of the
linear map $H\to\F$ by $h\mapsto h_a-h_b$.
So $r_{ij}=2^d$ or $2^{d-1}$ in this case.

Case: $i,j$ are of opposite same parity.  Let $i\in u_a$ and $j\in u_b$.
Subcase: $a=b$. Then $i,j$ are always either both in or both not in an $S_h$.
Then $S_h\oplus U$ will contain one of the $i,j$ but not the other.
This means $(a_i-a_j)$ will not appear in $\rho(\theta[\eta_{S_h}]^4)$
Thus $r_{ij}=0$ in this subcase.
Subcase: $a\ne b$. 
Then $\rho(\theta[\eta_{S_h}]^4)$ contains a factor of
$(a_i-a_j)$ if and only if $i,j$ are both in or both not in
$S_h\oplus U$ which happens if and only if exactly one of $i,j$
is in $S_h$ which is if and only if $h_a\ne h_b$.
Since the set of such $h$ is the complement in $H$
of the  kernel of the
linear map $H\to\F$ by $h\mapsto h_a-h_b$,
then $r_{ij}=|H|-|H|$ or
$|H|-\half |H|$.
So $r_{ij}=0$ or $2^{d-1}$ in this subcase.

In all cases, we
get $r_{ij}=0$, $2^{d-1}$, or $2^d$.
Note that $r_{ij}$ is a multiple of $2^{d-1}$ in all cases.
Thus there exists a polynomial
$Q_V$ such that $\rho(\prod_{\zeta\in V}\theta[\zeta]^4)=\pm
Q_V^{2^{d-1}}$.
Squaring both sides completes the proof.
\end{proof}

\begin{lm}\label{Q-lem}
Let $r,w\in\NN$.
 If $Q$ is a polynomial with real coefficients such that
$Q^4\in S_{4w}(r)$,
then $Q^2\in S_{2w}(r)$. 
\end{lm}

\begin{proof}
Take any $\smtwomat ABCD\in\SL_2(\RR)$ and $\gamma(z) =\frac{Az+B}{Cz+D}$.
We need to show that 
$Q(a_1,\ldots,a_r)^2 = Q(\gamma(a_1),\ldots,\gamma(a_r))^2\prod_{i=1}^r(Ca_i+D)^{2w}$
for all $a_i$.
We know that
\begin{equation}\label{Q-eq1}
Q(a_1,\ldots,a_r)^4 = Q(\gamma(a_1),\ldots,\gamma(a_r))^4\prod_{i=1}^r(Ca_i+D)^{4w}
\end{equation}
for all $a_i$.
Since $Q$ has real coefficients, then for real values of $a_i$,
Equation \ref{Q-eq1}  is an equality of real numbers to the fourth power.
Thus we can take the positive square root of both sides and obtain that
$Q(a_1,\ldots,a_r)^2 = Q(\gamma(a_1),\ldots,\gamma(a_r))^2\prod_{i=1}^r(Ca_i+D)^{2w}$
for all real values of the $a_i$.
That is, the multivariable polynomial
$Q(a_1,\ldots,a_r)^2 - Q(\gamma(a_1),\ldots,\gamma(a_r))^2\prod_{i=1}^r(Ca_i+D)^{2w}$
is zero for all real values of its variables and 
must be the zero polynomial.
Hence $Q^2\in S_{2w}(r)$.
\end{proof}

\begin{lm}\label{S-lem}
For any subspace $V$ of theta characteristics,
and any $\sigma\in\SS$, define
$$\sigma\cdot V
=\{\eta_{\sigma(S)}: \eta_S\in V \}.
$$
Then
$$\sigma(Q_V^2) = Q_{\sigma\cdot V}^2.$$
Furthermore, when $V$ is an isotropic subspace
consisting of balanced elements,
then $\sigma\cdot V$ is also such a subspace.
\end{lm}

\begin{proof}
For $\zeta=\eta_S\in V$, we have that
$\pm(a_i-a_j)$ occurs in $\rho(\theta[\zeta]^4)$
if and only if $i,j\in U_g\oplus S$ or $i,j\in U_g'\oplus S$.
Because $\sigma(U_g)=U_g$ or $\sigma(U_g)=U_g'$,
then $\pm(a_i-a_j)$ occurs in $\rho(\theta[\zeta]^4)$
if and only if $\sigma(i),\sigma(j)\in U_g\oplus \sigma(S)$ or 
$\sigma(i),\sigma(j)\in U_g'\oplus \sigma(S)$.
Thus 
$\pm(a_i-a_j)$ occurs in $\rho(\theta[\zeta]^4)$
if and only if
$\pm(a_{\sigma(i)}-a_{\sigma(j)})$ occurs in $\rho(\theta[\eta_{\sigma(S)}]^4)$.
Thus $\sigma(\rho(\theta[\zeta]^4))=\pm \rho(\theta[\eta_{\sigma(S)}]^4)$.
Letting $d=\dim V$,
we have by Lemma \ref{thelemma} that
$\sigma(Q_V^{2^d}) = Q_{\sigma\cdot V}^{2^d}$.
Since $Q_V$ and $Q_{\sigma\cdot V}$ have real coefficients, by an argument similar to that
of Lemma \ref{Q-lem},
we have that
$\sigma(Q_V^2) = Q_{\sigma\cdot V}^2.$

Now let $V$ be an isotropic subspace of balanced elements.
By Lemma \ref{balancelemma}, 
there exists a partitioning of $\{1,\ldots,2g+2\}$
into balanced subsets of 2 elements each (call
them $u_1,\ldots,u_{g+1}$),
and a subspace $H\subseteq \F^{g+1}$ with $\dim H=\dim V$
such that
$
V=\{\eta_{S_h}: h\in H, \text{ where }S_h=\bigcup_{i: h_i\ne0} u_i \}.
$
Then
$$
\sigma\cdot V=\{\eta_{\sigma(S_h)}: h\in H, \text{ where }
\sigma(S_h)=\bigcup_{i: h_i\ne0} \sigma(u_i) \},
$$
and so $\sigma\cdot V$ is a subspace,
and in fact a subspace of balanced elements.
\end{proof}

For any subspace $V$ of theta characteristics, we will
use the notation $Q_V$ as prescribed by Lemma \ref{thelemma},
with the understanding that $Q_V^2$ is unique given $V$.
Note that $Q_V=0$ unless $V$ is isotropic and contains only balanced elements.

Fix $g=g_1+g_2$ with $g_1,g_2\in\N$ for the following discussion, 
which parallels that of \cite{GR}.
For any theta characteristic $\zeta$,
the Witt map $W_{g_1,g_2}$ on Siegel modular forms yields
$$
W_{g_1,g_2} \theta[\zeta]^8
=
\theta[\pi_1\zeta]^8\theta[\pi_2\zeta]^8
$$
where
$\pi_1\zeta$ is the projection of $\zeta$ onto the left $2g_1$
coordinates
and
$\pi_2\zeta$ is the projection of $\zeta$ onto the right $2g_2$
coordinates.
Let $V\subseteq \F^g$ be a subspace of theta characteristics.
Then
$$
W_{g_1,g_2}\left( \prod_{\zeta\in V}\theta[\zeta]^8\right)
=
\prod_{\zeta_1\in \pi_1V}\theta[\zeta_1]^{8\cdot2^{d-d_1}}
\prod_{\zeta_2\in \pi_2V}\theta[\zeta_2]^{8\cdot2^{d-d_2}}
$$
where $d_i=\dim \pi_iV$.

Since the eighth powers are modular forms on $\Gamma(2)$
and the appropriate Witt maps,  $\Psi_{g_1,g_2}^*$ and $W_{g_1,g_2}$, 
are equivariant with respect to the $\rho$-map, 
we get that
$$
W_{g_1,g_2} \left( Q_V^4 \right)
=
Q_{\pi_1V}^{4\cdot2^{d-d_1}}
Q_{\pi_2V}^{4\cdot2^{d-d_2}}.
$$

Since $Q_V^4\in S_{8g}(2g+2)_0$
and $Q_V$ has real coefficients, 
then by Lemma \ref{Q-lem}, we have that
$Q_V^2\in S_{4g}(2g+2)_0$.  
The important point here is that $Q_V^2$ has the correct valuation
and we can apply the Witt map $W_{g_1,g_2}$ to it.
Then we must have
$$
W_{g_1,g_2} \left( Q_V^2 \right)
=
Q_{\pi_1V}^{2^{1+d-d_1}}
Q_{\pi_2V}^{2^{1+d-d_2}}.
$$

Now letting $V$ vary over subspaces of dimension d,
we get the following.
\begin{multline*}
W_{g_1,g_2} \left(\sum_{\substack{V\subseteq\F^{2g}\\\dim V=d}} Q_V^2 \right)
=
\sum_{0\le d_1,d_2\le d\le d_1+d_2}
N_{d_1,d_2;d}\\
\cdot
\sum_{\substack{V_1\subseteq\F^{2g_1}\\\dim V_1=d_1}}
Q_{V_1}^{2^{1+d-d_1}}
\sum_{\substack{V_2\subseteq\F^{2g_2}\\\dim V_2=d_2}}
Q_{V_2}^{2^{1+d-d_2}},
\end{multline*}
where $N_{d_1,d_2;d}$ is the number of $V\subseteq\F^{g_1}\oplus\F^{g_2}$ 
of dimension $d$ that have $\pi_iV=V_i$ ($i=1,2$) given
fixed $V_1,V_2$ of dimensions $d_1,d_2$ respectively.
The formula proven in \cite{GR} is
$$N_{d_1,d_2;d}=\prod_{j=0}^{d_1+d_2-d-1}\frac{(2^{d_1}-2^j)(2^{d_2}-2^j)}{(2^{d_1+d_2-d}-2^j)}$$
for $0\le d_1,d_2\le d\le d_1+d_2$ and it is $0$ otherwise.
Now define
$$
\tilde G^{(g)}_{d}=
\sum_{\substack{V\subseteq\F^{2g}\\\dim V=d}} Q_V^2.
$$
Note that by summing over all subspaces, by Lemma \ref{S-lem}, we have
that $\tilde G^{(g)}_{d}$ is invariant under $\SS$.  
We have
\begin{equation}\label{WGG}
W_{g_1,g_2}\tilde G^{(g)}_d =
\sum_{0\le d_1,d_2\le d\le d_1+d_2}
N_{d_1,d_2;d}\,
\tilde G^{(g_1)}_{d} \tilde G^{(g_2)}_{d}
\end{equation}

Define
\begin{equation}\label{Kg}
K^{(g)}=
\frac1{2^g}\left(\sum_{i=0}^g(-1)^i 2^{\frac{i(i-1)}2}\tilde G^{(g)}_{i}
\right).
\end{equation}

\begin{prop}
Let $g=g_1+g_2$ with $g_1,g_2\in\N$.
Then
$$W_{g_1,g_2}K^{(g)} = K^{(g_1)}  \otimes K^{(g_2)}$$
\end{prop}

\begin{proof}  
Expand both sides using Equations \ref{WGG} and \ref{Kg}
 and prove that the coefficient of 
$\tilde G^{(g_1)}_{n}
\tilde G^{(g_2)}_{m}$
is the same for all $n,m$.
The key identity is
$$
(-1)^n 2^{\frac{n(n-1)}2}
(-1)^m 2^{\frac{m(m-1)}2}
=
\sum_{i=0}^{n+m}(-1)^i 2^{\frac{i(i-1)}2}N_{n,m;i}
$$
which is proven in \cite{GR}.
\end{proof}

\begin{thm}\label{altthm}
The family
$K^{(g)}$ satisfies the hyperelliptic Ansatz.
Furthermore,
$K^{(g)}=H_g$ for all $g$.
\end{thm}

\begin{proof}
We already know that $K^{(g)} \in S_{4g}(2g+2)_0(\SS)=\BB_g^8$ is a family of
modular forms that satisfy the splitting property
and that $K^{(1)}$ satisfies the base condition.
By Theorem \ref{unique} on uniqueness, we must have
that the family $K^{(g)}$ equals the family $H_{g}$.
\end{proof}

\end{document}